\documentclass[12pt]{article}
\usepackage{amsmath}
\usepackage{amssymb}
\usepackage{cite}

\newcommand{\eT}{{\mathbf{T}}}

\def\nn{\nonumber}

\def\nn{\nonumber}

\def\b{\beta}
\def\g{\gamma}

\def\vk{\varkappa}

\def\ts{\times}

\def\vp{\varphi}
\def\ve{\varepsilon}
\def\wh{\widehat}
\def\wt{\widetilde}
\def\ov{\overline}
\def\mg{{\mathfrak G}}
\def\BC{{\mathbb C}}
\def\BR{{\mathbb R}}
\def\clp{{\mathcal P}}

\def\cln{{\mathcal N}}
\def\clu{{\mathcal U}}

\def\const{{\rm const}}

\newtheorem{Pa}{Paper}[section]
\newtheorem{Tm}[Pa]{{\bf Theorem}}
\newtheorem{La}[Pa]{{\bf Lemma}}
\newtheorem{Cy}[Pa]{{\bf Corollary}}
\newtheorem{Rk}[Pa]{{\bf Remark}}

\newtheorem{Dn}[Pa]{{\bf Definition}}
\newtheorem{Pn}[Pa]{{\bf Proposition}}

\newcommand{\iy}{{\infty}}

\newcommand{\bpr}{{\noindent\textbf{Proof.}\ }}
\newcommand{\epr}{{\hfill $\Box$}}

\newcommand{\E}{\mathrm{e}}
\newcommand{\I}{\mathrm{i}}

\title{Inverse problem for Dirac systems with locally square-summable potentials and rectangular Weyl
functions}

\author{Alexander Sakhnovich}

\date{}

\begin{document}

\maketitle

\thanks{
Address: Fakult\"at f\"ur Mathematik, Universit\"at Wien,\\
Oskar-Morgenstern-Platz 1, A-1090 Wien, Austria;\\
e-mail:  oleksandr.sakhnovych@univie.ac.at}

\begin{abstract} Inverse problem for Dirac systems with locally
square summable potentials and rectangular Weyl
functions is solved. For that purpose we use a new result
 on the linear similarity between operators from a subclass of triangular integral
 operators and the operator of integration.
\end{abstract}



\section{Introduction}
We consider
the self-adjoint Dirac (more precisely, Dirac-type) 
system
\begin{align} &       \label{1.1}
\frac{d}{dx}y(x, z )=\I (z j+jV(x))y(x,
z ) \quad
(x \geq 0),
\end{align} 
where
\begin{align} &   \label{1.2}
j = \left[
\begin{array}{cc}
I_{m_1} & 0 \\ 0 & -I_{m_2}
\end{array}
\right], \hspace{1em} V= \left[\begin{array}{cc}
0&v\\v^{*}&0\end{array}\right],  \quad m_1+m_2=:m,
 \end{align} 
$I_{m_k}$ is the $m_k \times
m_k$ identity
matrix and $v(x)$ is an $m_1 \times m_2$ matrix function.
 We  assume
that  $v$ is measurable and, moreover, locally square-summable, that is, square-summable
on  the finite intervals $[0, \, l]$. Here we say that a matrix function is summable (square-summable) if its entries are summable (square-summable). 

Dirac (Dirac-type) system is a classical object of analysis. Its Weyl and spectral theories were actively studied in the second half of the 20-th
century, the first solution of the inverse spectral problem being given (for the case of the scalar $v$ and without proof) 
by M.G. Krein in the seminal paper \cite{Kre1}. For the quite recent publications on Dirac systems see, for instance, 
\cite{AHM, AGKLS2, CG, FKRS3, FKRSp1, GGHT, MyPuy, Puy, SaA3, SaSaR}
and references therein. Dirac system is of independent interest and it is also important as an auxiliary system for many integrable nonlinear
equations. Moreover, it is related to the famous Schr\"odinger equation (see, e.g., \cite{EGNT}).  Many recent publications are dedicated to
the development of the Weyl and spectral theories of Dirac system under weaker summability conditions. Here, we solve the inverse problem
under the  condition of the   local square-summability of $v$. We deal with the case, where the potential  $v$ and the corresponding
Weyl function are rectangular (not necessarily square)  matrix functions, which is essential for some applications to the matrix and
multicomponent integrable equations.

Before stating our main result, we formulate several  results from \cite{FKRSp1, SaSaR} on direct problems.
The notation $u(x,z)$ stands for the fundamental solution of \eqref{1.1}
normalized by the condition
\begin{align} &   \label{1.3}
u(0,z)=I_m.
\end{align}
Later we shall need notations of the block rows of $u(x,0)$:
\begin{align} &      \label{3.1}
\b(x)=\begin{bmatrix}
I_{m_1} & 0
\end{bmatrix}u(x,0), \quad \g(x)=\begin{bmatrix}
0 &I_{m_2}
\end{bmatrix}u(x,0).
\end{align} 
\begin{Dn} \label{defWeyl} Weyl-Titchmarsh $($or simply Weyl$)$ function of Dirac system  \eqref{1.1} on $[0, \, \infty)$,
where  the potential $v$ is locally summable, is a holomorphic $m_2\times~m_1$ matrix function $\vp$
which satisfies the inequality
\begin{align} &      \label{2.20}
\int_0^{\infty}
\begin{bmatrix}
I_{m_1} & \vp(z)^*
\end{bmatrix}
u(x,z)^*u(x,z)
\begin{bmatrix}
I_{m_1} \\ \vp(z)
\end{bmatrix}dx< \infty , \quad z\in \BC_+.
\end{align} 
\end{Dn}
Here $\BC_+$ stands for the open upper half-plane.
In order to study Weyl functions, we introduce the class  of nonsingular $m \times m_1$ matrix functions 
$\clp(z)$ with property-$j$. Namely, the matrix functions 
$\clp(z)$ are meromorphic in $\mathbb{C}_+$ and satisfy
(excluding, possibly, a discrete set of points)
the following relations
\begin{align}\label{2.1}&
\clp(z)^*\clp(z) >0, \quad \clp(z)^* j \clp(z) \geq 0 \quad (z\in \mathbb{C}_+).
\end{align} 
Relations  (\ref{2.1})  imply
\begin{align} &      \label{2.7}
\det \Big(\begin{bmatrix}
I_{m_1} & 0
\end{bmatrix}u(x,z)^{-1}\clp(z)\Big)\not= 0.
\end{align} 
\begin{Dn} \label{set}
The set $\cln(x,z)$ of M\"obius transformations is the set of values at $x, \,z$ 
of matrix functions
\begin{align}\label{2.2}&
\vp(x,z,\clp)=\begin{bmatrix}
0 &I_{m_2}
\end{bmatrix}u(x,z)^{-1}\clp(z)\Big(\begin{bmatrix}
I_{m_1} & 0
\end{bmatrix}u(x,z)^{-1}\clp(z)\Big)^{-1},
\end{align} 
where $\clp(z)$ are nonsingular  matrix functions 
 with property-$j$. 
 \end{Dn}
 As usual, the sets $\cln(x,z)$ are embedded, that is, 
 \begin{align}\label{2.2'}&
 \cln(x_1,z)\subset \cln(x_2,z) \quad {\mathrm{for}}\quad x_1>x_2.
 \end{align}
 Moreover, the following proposition holds.
\begin{Pn} \label{CyW2} \cite[Subsection 2.2.1]{SaSaR}
Let Dirac system  \eqref{1.1} be given on $[0, \, \infty)$
and let its  potential $v$ be locally summable. Then
there is a unique matrix function $\vp(z)$ in $\BC_+$ such that
\begin{align} &      \label{2.3}
\vp(z)=\bigcap_{x<\infty}\cln(x,z).
\end{align}  
This function is analytic and non-expansive. Moreover, this function is the unique
Weyl function of system \eqref{1.1}.
\end{Pn}
If $v$ is locally square-summable, we may recover it from the Weyl function.
\begin{Tm}\label{TmMai}
Let Dirac system  \eqref{1.1} be given on $[0, \, \infty)$,
let its  potential $v$ be locally square-summable and let
$\vp$ be the Weyl function of this system. Then $v$ is uniquely recovered
from $\vp$.
\end{Tm}
The procedure to recover $v$ from $\vp$ is based on the study of the operator
\begin{align} &      \label{3.4}
K=\I \int_0^x  \, \g(x)j\g(t)^* \cdot \, dt, \quad K \in B\big(L^2_{m_2}(0, \, l)\big),
\end{align} 
where $\g$ is the lower block row of $u(x,0)$ (see \eqref{3.1})
and $B(H)$ denotes the class of bounded
linear operators, which map the space $H$ into $H$.
Using a new version of the similarity result for $K$, we modify
the procedure to solve inverse problem, which was  developed in \cite{SaA02, SaA3, SaSaR}, for the case
of the less smooth than before potentials $v$.

Further $F^{\prime}$ stands for the derivative of $F$, "$\const$" means a constant function or vector-function,
$I_r$ is the $r\times r$ identity matrix, $I$ is an identity operator, $B(H_1, H_2)$ denotes the class of bounded
linear operators, which map the Hilbert space $H_1$ into the Hilbert space $H_2$.
Speaking about fundamental solutions we assume that they are normalized by $I_m$ at $x=0$.


\section{Similarity result}
\setcounter{equation}{0}
We consider conditions of similarity of the two operators acting in $L^2_r(0, \,\eT)$, namely,
\begin{equation} \label{a10}
K:=F(x)\int_0^x G(t)\, \cdot \,dt, \quad A:=\int_0^x \, \cdot \,dt,
\end{equation}
where $F$ and $G$ are differentiable $r \times p$ and $p \times r$, respectively, matrix functions.
\begin{Pn} \label{PnSim1} Let $F$ and $G$ be differentiable and satisfy the identity
\begin{align}\label{n1}
F(x)G(x) \equiv I_r, \quad 0 \leq x \leq {\mathbf T},
\end{align}
and assume that the entries of $F^{\prime}$ and $G^{\prime}$ belong $L^2(0, \, \eT)$.

Then the operator $K$ defined by \eqref{a10}
is similar to the operator of integration $A$. More precisely, $K=E A
E^{-1}$ where $E \in B\big(L^2_r(0, \, \eT)\big)$ is a lower triangular operator of the form
\begin{equation}
\label{a27} E= \rho(x)\left(I+\int_0^x N(x,t)\, \cdot \,dt\right), \quad \frac{d}{dx}\rho=F'G\rho, \quad \rho(0)=I_r,
\end{equation}
and the matrix functions $\rho$, $\rho^{-1}$ and $N$ are measurable and uniformly bounded.
Moreover, the operators $E^{\pm 1}$ map differentiable functions with a square-summable derivative
into differentiable functions with a square-summable derivative.
\end{Pn}
The case of operators $K$ of the form \eqref{a10}, where $F$ and $G$ have bounded derivatives, is a particular
case of operators, the similarity of which to $A$ was proved in an important paper \cite{SaL00}.
Later on, the proof from \cite{SaL00} was modified for the case of operators $K$ such that $F$ and $G$
have continuous derivatives (and $E^{\pm 1}$ map functions with continuous derivatives into functions
with continuous derivatives) \cite{AGKLS2}. Here, we modify further the proofs from \cite{AGKLS2, SaL00}
for the case of the less smooth functions $F$ and $G$.
The proof of  Proposition \ref{PnSim1} above requires some preparations. 

We note that, according to the general
theory of  semi-separable integral
operators, which is also easily checked directly, the inverse of operator $I - z K$
is given by
\begin{equation} \label{na15}
\big((I - z K)^{-1}f\big)(x)=
f(x)+\int_0^x Q(x,t,z)  f(t) dt,
\end{equation}
where  
\begin{eqnarray} \label{na16}&&
  Q(x,t,z )=zF(x)u_1(x,z)u_1(t,z)^{-1}G(t), \quad 0 \leq t \leq x \leq {\bf T}; \\
  \label{1b2}
  &&  \frac{d}{d x}\, u_1(x,z)= zG(x)F(x)u_1(x,z), \quad  0 \leq x \leq {\bf T}; \\
\label{1b3}
&&  u_1(0,z)=I_r.
\end{eqnarray}
Introduce also the $p \times p$ matrix function $\wt u_1(x)$ defined by
\begin{equation} \label{1b1}
\frac{d}{d x}\wt  u_1(x)= - G(x)F^{\prime}(x) \wt u_1(x), \quad  0 \leq x \leq {\bf T}, \\
\quad \wt u_1(0)=I_p.
\end{equation}
We are now ready to prove the first lemma.

\begin{La} \label{inteq}
Let $F$  and  $G$ be absolutely continuous and assume that the identity
\eqref{n1}
holds.
Introduce the $r\times r$ matrix functions $h$ and $\rho$ by 
\begin{align}\label{defrho}
h(x):=F(x)G(0); \quad \frac{d}{dx}\rho=F'G\rho, \quad \rho(0)=I_r.
\end{align}
Put
\begin{equation} \label{a13}
g(x, z)=\rho(x)^{-1}\left((I-zK)^{-1}h\right)(x), \quad 0\leq x\leq \eT,
\end{equation}
where $(I-zK)^{-1}$ is applied to $h$ columnwise. Then $g$ satisfies the
following integro-differential equation
\begin{equation} \label{a20}
\frac{d}{d x}\,g(x, z)-\mu(x)\int_0^x \nu(t) g(t, z)dt -zg(x, z)=0, \quad g(0, z)=I_r,
\end{equation}
where $\mu$ and $\nu$  are the summable functions on  $[0,\,\eT]$ given by
\begin{align}
\mu(x):&=\rho(x)^{-1}F^{\prime}(x)\wt u_1(x),\quad  0 \leq x \leq \eT; \label{defagr}\\
\nu(t):&=-\wt u_1(t)^{-1}\big(G(t)F^{\prime}(t)G(t)+G^{\prime}(t)\big)\rho(t),\quad 0 \leq t \leq \eT.\label{defbgr}
\end{align}
\end{La}
\bpr
Put $\wt g(x, z)=\rho(x)g(x, z)$. Using (\ref{na15})-(\ref{1b3}), (\ref{a13}), and the definition of the matrix
function $h$, we present $\wt g$ in the form
\begin{align}\nn
\wt g(x, z)&=  F(x)G(0)
\\ & \nn \,\,\quad+  zF(x)u_1(x,z)\int_0^x u_1(t,z)^{-1}G(t) F(t)G(0)dt  \\
&=F(x)G(0)-F(x)u_1(x,z)\int_0^x \frac{d}{dt}\Big(u_1(t,z)^{-1}G(0)\Big)dt\nonumber \\
& = F(x)G(0)-F(x)u_1(x,z)\big(u_1(x,z)^{-1}-I_r\big)G(0)\nonumber\\
&=F(x)u_1(x,z)G(0). \label{4a1}
\end{align}
It follows that
\begin{equation} \label{4a2}
g(x, z)=\rho(x)^{-1}F(x)u_1(x,z)G(0).
\end{equation}
Clearly  $g$ is differentiable and
\begin{eqnarray} \label{4a3}&&
\frac{d}{d x}\,g(x, z)=\rho(x)^{-1}\wt g_x(x, z)-  \rho(x)^{-1}\rho^{\prime}(x)\rho(x)^{-1}  \wt g(x, z)
\\
\noalign{\vskip4pt}
&&=\rho(x)^{-1} \big\{
zF(x)G(x)F(x)+F^{\prime}(x) -F^{\prime}(x) G(x)F(x)\big\}u_1(x, z)G(0)
    \nonumber \\
\noalign{\vskip4pt}
  &&=zg(x, z)+\rho(x)^{-1} F^{\prime}(x) \big(I_p-G(x)F(x)\big)u_1(x, z)G(0). \nonumber
\end{eqnarray}
Here we took into account the identity \eqref{n1}.
From \eqref{1b1} we see that
\[
\frac{d}{dt}  \wt u_1(t)^{-1}=- \wt u_1(t)^{-1}\left(\frac{d}{dt} \wt u_1(t)\right)\wt u_1(t)^{-1}=\wt u_1(t)^{-1}G(t)F^\prime(t).
\]
Hence
\begin{eqnarray*}
&&\frac{d}{dt}  \left(\wt u_1(t)^{-1}\big(I_p-G(t)F(t)\big)
u_1(t,z)\right)\\
&&\hspace{1.5cm}  =\wt u_1(t)^{-1}G(t)F^{\prime}(t)\big(I_p-G(t)F(t)\big)u_1(t,z)\\
&&\hspace{2.5cm}+\wt u_1(t)^{-1}\big(-G^{\prime}(t)F(t)- G(t)F^{\prime}(t)\big)u_1(t,z)\\
&&\hspace{2.5cm}+z\wt u_1(t)^{-1}\big(I_p-G(t)F(t)\big)G(t)F(t)u_1(t,z).
\end{eqnarray*}
Since, in view of condition \eqref{n1}, we have $\big(I_p-G(t)F(t)\big)G(t)=0$, we obtain 
\begin{eqnarray*}
&&\frac{d}{dt}  \left(\wt u_1(t)^{-1}\big(I_p-G(t)F(t)\big)
u_1(t,z)\right)\\
&&\hspace{1.5cm} =\wt u_1(t)^{-1}\big(G(t)F^{\prime}(t)-G(t)F^{\prime}(t)G(t)F(t)\\
&&\hspace{3.5cm}-G^{\prime}(t)F(t)- G(t)F^{\prime}(t)\big)u_1(t,z)\\
&&\hspace{1.5cm}=-\wt u_1(t)^{-1}\big(G(t)F^{\prime}(t)G(t)+G^{\prime}(t)\big)F(t)u_1(t,z).
\end{eqnarray*}
Using the definition of $\nu$  in \eqref{defbgr} and the identity \eqref{4a2}, we derive
\begin{equation}\label{sum1}
\frac{d}{dt}  \left(\wt u_1(t)^{-1}\big(I_p-G(t)F(t)\big)
u_1(t,z)\right)G(0)=\nu(t)g(t,z).
\end{equation}
Recall that  $\big(I_p-G(t)F(t)\big)G(t)=0$ and so $\big(I_p-G(0)F(0)\big)G(0)=0$, in particular.
Hence, from \eqref{sum1} it follows that
\begin{align}\nn
\int_0^x \nu(t)g(t,z)\,dt&=\wt u_1(x)^{-1}\big(I_p-G(x)F(x)\big)u_1(x,z)G(0)
-\big(I_p-G(0)F(0)\big)
\\ \label{sum2}
\times G(0)
&
= \wt u_1(x)^{-1}\big(I_p-G(x)F(x)\big)u_1(x,z)G(0).
\end{align}
But then, using \eqref{4a3} and the definition of $\mu$ in \eqref{defagr}, we arrive at the
identity \eqref{a20}.
\epr

The lemma below provides an integral representation  of the solution of \eqref{a20}.
\begin{La} \label{repr}
Let $\mu(x)$ and $\nu(x)$ be  $r \times p$ and $p \times r$, respectively, matrix functions,
such that  their entries  belong $L^2(0, \,\eT)$.
Then the
integro-differential equation
\begin{align}& \label{b20}
\frac{d}{d x}\,g(x, z)-\int_0^x \vk(x,t) g(t, z)dt -zg(x, z)=0, \quad g(0, z)=I_r, \\
&\label{b20'} \vk(x,t):=\mu(x)\nu(t)
\end{align}
has a unique
solution $g(\cdot, z)\in L^2_r(0, \, \eT)$, and this solution has the form
\begin{equation}
\label{n3}g(x,z)=e^{zx}I_r+\int_0^x e^{zt} N(x,t) \,dt,\quad
0\leq x\leq \eT,
\end{equation}
where $N(x,t)$ is  bounded on $0\leq t\leq x\leq\eT$. 
\end{La}
\bpr  We set
\begin{align}
\label{a24}&
\vk_1(x,t)=\int_{x-t}^x\vk(\xi,\xi+t-x)d\xi, \quad 0\leq t\leq x\leq\eT,
\\ \label{a25}&
\vk_{k+1}(x,t)=\int_{x-t}^x
\int_{y+t-x}^y \vk(y,s)\vk_k(s,y+t-x)\,ds\,dy.
\end{align}
It is easily proved by induction that
\begin{equation}\label{gk2}
\|\vk_k(x,t)\|\leq C_0C_1^{k-1}\frac{x^{k-1}}{(k-1)!}, \quad 0\leq t \leq
x\leq \eT, \quad k\geq 1
\end{equation}
for some $C_0, C_1>0$. Thus, we can introduce a bounded matrix function
\begin{align}
\label{sum3}&
N(x,t)=\sum_{k=1}^\iy \vk_k(x,t), \quad 0\leq t\leq x\leq \eT.
\end{align}
Putting 
\begin{align}
\label{sum4} &
\mg_0(x,t)=\vk(x, t); \quad \mg_k(x, t)=\int_{t}^x \vk(x,s)\vk_k(s,t)\,ds, 
 \quad k>0,
\end{align}
and using \eqref{a24}, \eqref{a25},  and \eqref{sum4}, we easily derive
\begin{align}
\label{sum5}
\int_0^x e^{z(x-\xi)}\left(\int_0^\xi e^{zt}\mg_k(\xi, t)dt\right)  d\xi &=
\int_0^x \left(\int_0^\xi e^{z(x+t-\xi)}\mg_k(\xi,t)\, dt\right)  d\xi \\& \nn
=\int_0^x\left(\int_{x-\xi}^x e^{zt}\mg_k(\xi,\xi+t-x) dt\right)  d\xi\\
&\nn =\int_0^x e^{zt}\left(\int_{x-t}^x\mg_k(\xi,\xi+t-x)d\xi\right)dt
\\& \nn
= \int_0^x e^{zt}\vk_{k+1}(x,t)\, dt \quad (k\geq 0).
\end{align}
Taking into account \eqref{sum3} and \eqref{sum5}, we see that $g$ given by \eqref{n3} satisfies
the equation
\begin{align}& \label{y1}
\frac{d}{d x}\,g(x, z) -zg(x, z)=\int_0^x e^{zt}\left(\sum_{k=0}^{\infty}\mg_k (x,t) \right)dt.
\end{align}
In view of \eqref{sum4} we have the equalities
\begin{align} \label{y2} &
\int_0^x e^{zt}\vk (x,t) dt=\int_0^x e^{zt}\mg_0 (x,t) dt, \\ \nn &
\int_0^x\vk(x,t)\int_0^t e^{zs}\vk_k(t,s)dsdt=\int_0^x\vk(x,s)\int_0^s e^{zt}\vk_k(s,t)dtds
\\  \label{y3}&
=\int_0^x e^{zt} \int_t^x\vk(x,s) \vk_k(s,t)dsdt
= \int_0^x e^{zt}\mg_k (x,t) dt.
\end{align}
Using \eqref{n3}, \eqref{sum3}, \eqref{y2}, and \eqref{y3}, we rewrite \eqref{y1}
in the form \eqref{b20}. 

It remains to prove that the solution of \eqref{b20} is unique. Indeed, integrating \eqref{b20}
with respect to $x$ we derive the equality
\begin{align} \label{y4} &
g(\cdot, z)-ARg(\cdot, z)-zAg(\cdot, z)=I_r,
\end{align}
where the bounded in $L^2_r(0, \eT)$ operators $A$ and $R$ are given by the relations
\begin{align} \label{y5} &
Af=\int_0^xf(t)dt, \quad Rf=\int_0^x\vk(x,t) f(t)dt.
\end{align}
Clearly $A$ is a Volterra operator and it is easily checked (see also, e.g., \cite[Subsection 1.2.4]{SaSaR} and \cite{SaLbez})
that 
\begin{align} \label{y5'} &
(I-zA)^{-1}= I+z\int_0^xe^{z(x-t)}\cdot dt.
\end{align}
Therefore, $(I-zA)^{-1}AR$ is an integral triangular operator
with Hilbert-Schmidt kernel (and so $(I-zA)^{-1}AR$ is also a Volterra operator).
Hence, according to \eqref{y4}, the solution $g$ of \eqref{b20} is uniquely
defined by the formula
\begin{align} \label{y6} &
g(\cdot, z)=(I-(I-zA)^{-1}AR)^{-1}(I-zA)^{-1}I_r.
\end{align}
\epr

{\it Proof of Proposition \ref{PnSim1}}. We split the proof into two steps. In the first step
we construct the operator $E$ and establish the similarity $KE=EA$. In the next step we prove that
$E^{\pm 1}$ map
functions with a square-summable derivative
into functions with a square-summable derivative.

Step 1. Let $g(x,z)$ be the matrix function defined by \eqref{a13}.
According to Lemma \ref{inteq},  $g(x,z)$ satisfies the equation
\eqref{a20}. Hence, in view of Lemma \ref{repr}, 
$g$ admits the representation
\begin{equation}
\label{n3c}
g(x,z)=e^{zx}I_r+\int_0^x N(x,t)\big(e^{zt}I_r\big)dt,\quad
0\leq x\leq \eT,
\end{equation}
where $N(x,t)$ is given by \eqref{sum3}. The same $N(x,t)$ is substituted into the definition
\eqref{a27} of the operator $E$ acting on  $L^2_r(0, \,\eT)$, whereas
the $r\ts r$ matrix function $\rho$ in \eqref{a27} coincides with $\rho$ defined by \eqref{defrho}. 
Thus, the matrix functions $\rho$, $\rho^{-1}$ and $N$ are measurable and uniformly bounded,
and $E$ is boundedly invertible.

Taking into account \eqref{a27},  \eqref{a13}, and \eqref{n3c}
we see that 
\begin{align} \label{y7} &
E\big(e^{z x}I_r\big)=\rho(x)g(x,z)=(I-zK)^{-1}h,
\end{align}
where $h$ is determined in \eqref{defrho} (i.e., $h(x)=F(x)G(0)$).
It is immediate from \eqref{y5'} that
\begin{align} \label{y8} &
e^{z x}I_r=(I-zA)^{-1}I_r.
\end{align}
For the case that $z=0$ formula \eqref{y7} yields $EI_r=h$.
Thus, using \eqref{y8}, we rewrite \eqref{y7} in the form
\begin{align} \label{y9} &
E(I-zA)^{-1}I_r=(I-zK)^{-1}E I_r.
\end{align}
From the series expansion in (\ref{y9}) it follows that
\begin{equation} \label{a31}
E A^{j}I_r=K^j E I_r, \qquad j=0,1,2, \ldots.
\end{equation}
Therefore, for each $j=0,1,2,\ldots$, we have
\begin{equation} \label{a32}
(K E)A^{j}I_r=K(EA^{j}I_r)=K^{j+1}E I_r=E A^{j+1}I_r=(E A )A^{j}I_r.
\end{equation}
As the closed linear span  of the columns of the matrices
$\{A^jI_r\}_{j=0}^{\infty}$ coincides
with $L^2_r(0,\, \eT)$, the equalities  in (\ref{a32}) yield  $KE=EA$. Since $E$ is invertible, we obtain
$K=EAE^{-1}$, and hence  $K$ and $A$ are similar. It remains to prove that $E^{\pm 1}$ map
functions with a square-summable derivative
into functions with a square-summable derivative.

Step 2. Let $f$ be a differentiable vector function such that \\ $\wt f:=f^{\prime}\in L^2_r(0, \, \eT)$.
Then $f$ admits a representation
\begin{align} \label{y10} &
f=A\wt f +f_0 \quad (\wt f\in L^2_r(0, \, \eT)), \quad f_0\equiv \const .
\end{align}
According to the previous step,  $EI_r=h(x)=F(x)G(0)$, and so
\begin{align} \label{y10'} &
Ef_0=F(x)G(0)f_0, \quad  (Ef_0)^{\prime}=F^{\prime}(x)G(0)f_0.
\end{align}
Since we assume that the derivative $F^{\prime}$ is square-summable,
the same is valid for $Ef_0$. Next note that
\begin{align} \label{y11} &
(EA \wt f)(x)=(KE\wt f)(x)=F(x)\int_0^x G(t)(E\wt f)(t)\,dt.
\end{align}
Since 
$E$ maps $L^2_r(0, \, \eT)$ onto $L^2_r(0, \, \eT)$, 
 formula \eqref{y11} shows that $EA\wt f$  has a square-summable derivative.
 Thus, both $Ef_0$ and $EA\wt f$   have square-summable derivatives.
 Therefore, \eqref{y10} implies that $Ef$ also has a square-summable derivative.

Finally, we consider $E^{-1}$. First, introduce operator 
$K_1$  on $L^2_r(0,\,\eT)$:
\[
(K_1f)(x)=F^\prime(x)\int_0^xG(t)f(t)dt, \quad f\in L^2_r(0,\,\eT),
\]
and notice that $AK_1=K-A$ or, equivalently,
\begin{align} \label{y12} &
A(I+K_1)=K.
\end{align}
The operator $K_1$ is a triangular operator with  Hilbert-Schmidt kernel. In particular,
$K_1$ is a Volterra operator. Thus, $I+K_1$ is invertible.
Since $E$ is also invertible, we rewrite
$KE=EA$ as $E^{-1}K=AE^{-1}$. In view of \eqref{y12} the equality  $E^{-1}K=AE^{-1}$ yields
\begin{eqnarray}
E^{-1}A&=&E^{-1}A(I+K_1)(I+K_1)^{-1}=E^{-1}K(I+K_1)^{-1}\label{a35}\\
&=&AE^{-1}(I+K_1)^{-1}.\nonumber
\end{eqnarray}

Recall that $f$ with a square-summable derivative admits the representation \eqref{y10}.
Formula \eqref{a35} implies that $E^{-1}A\wt f$   has a square-summable derivative.
In order to show that $E^{-1} f_0$  also has a square-summable derivative, we 
take into account \eqref{n1} and rewrite the first equality
in \eqref{y10'} in the form
$$f_0=E^{-1}\big(F(x)G(0)f_0\big)=E^{-1}A\big(F^{\prime}(x)G(0)f_0\big)+E^{-1}f_0,$$
that is, 
\begin{align} \label{y13} &
E^{-1}f_0=f_0-E^{-1}A\big(F^{\prime}(x)G(0)f_0\big),
\end{align}
which completes the proof.

\epr
\begin{Rk}\label{Rk0} Relations \eqref{y10}, \eqref{a35}, and \eqref{y13} show that for any differentiable $f$
with a square-summable derivative
we have
\begin{align} \label{i8} &
(E^{-1}f)(0)=f(0).
\end{align}

\end{Rk}

\section{Dirac system: fundamental solution}
\setcounter{equation}{0}
We start with a similarity result, which follows from Proposition \ref{PnSim1}.
\begin{Pn} \label{PnSimN} Let  the potential $v$ of Dirac system \eqref{1.1} be 
square-summable on $(0, \, \eT)$, and let
$K$ be given by  \eqref{3.4}, where $\g$ is defined in \eqref{3.1}.
Then there is a similarity transformation operator $E \in B\big(L^2_r(0, \, \eT)\big)$ such that
\begin{align}\label{i1}
& K=EAE^{-1}, \quad  A:=-\I \int_0^x \, \cdot \,dt,
\\ & \label{i2} E= I+\int_0^x N(x,t)\, \cdot \,dt, \quad
\\  & \label{i3} E^{-1}\g_2\equiv I_{m_2},
\end{align}
where 
 $N$ is a Hilbert-Schmidt kernel
 and $\g_2$ is the right $m_2\times m_2$ block of $\g$.
Moreover, the operators $E^{\pm 1}$ map differentiable functions with a square-summable derivative
into differentiable functions with a square-summable derivative.
\end{Pn}
\bpr 
According to \eqref{1.1} we have
\begin{align}\label{i4}
u(x,0)^*ju(x,0)=j=u(x,0)ju(x,0)^*. 
\end{align}
Therefore, the blocks of $u(x,0)$ introduced in \eqref{3.1} satisfy the relations
\begin{align} &      \label{i5}
\b j \b^*\equiv I_{m_1}, \quad \g j \g^*\equiv -I_{m_2}, \quad \b j \g^*\equiv 0. 
\end{align} 
Furthermore, equation \eqref{1.1} implies that $\g^{\prime}$ is square-summable and
$$\g^{\prime}(x)=-\I \begin{bmatrix}v(x)^* & 0\end{bmatrix}u(x,0)=-\I v(x)^*\b(x).$$
Hence, the third equality in \eqref{i5} yields
\begin{align} &      \label{i6}
 \g^{\prime} j \g^*\equiv 0.
\end{align} 
In view of the second equality in \eqref{i5}, we may apply Proposition \ref{PnSim1} to $\I K$
(where $K$ is defined in \eqref{3.4}).
Moreover, \eqref{i6} implies the indentity $\rho(x) \equiv I_r$ for $\rho$ given in \eqref{a27}.
Thus, there is some similarity transformation operator $\wt E$, which satisfies all
conditions of Proposition \ref{PnSimN} excluding, possibly, equality \eqref{i3}
(and the kernel of $\wt E$ is bounded).
Let us normalize $\wt E$ multiplying it by the operator
\begin{align} &      \label{i7}
E_0=I + \int_0^x E_0(x-t) \cdot dt, \quad 
E_0(x):=\big(\wt E^{-1}\g_2\big)^{\prime}(x).
\end{align} 
We see that $E=\wt E E_0$ admits representation \eqref{i2}, where $N$
is a Hilbert-Schmidt kernel and that $AE_0=E_0A$. Thus, from $K=\wt EA \wt E^{-1}$
follows $K= EA  E^{-1}$. Finally, in view of \eqref{i7} and Remark \ref{Rk0} we obtain
\begin{align}    \nn
\big(E_0 I_{m_2}\big)(x)&=
I_{m_2}+\int_0^xE_0(t)dt
=I_{m_2}+
\big(\wt E^{-1} \g_2 \big)(x)-\big(\wt E^{-1} \g_2 \big)(0)
\\ \label{p3}&
=\big(\wt E^{-1} \g_2 \big)(x),
\end{align} 
and so \eqref{i3} is valid for $E=\wt E E_0$.

Clearly, the equalities $AE_0=E_0A$ and \eqref{p3} imply
that $E_0$ maps differentiable functions with a square-summable derivative
into differentiable functions with a square-summable derivative.
Rewriting $AE_0=E_0A$ and \eqref{p3} in the forms
\begin{align}    \nn
E_0^{-1}A=AE_0^{-1}, \quad E_0^{-1}I_{m_2}=I_{m_2}-\I E_0^{-1} A \big(\wt E^{-1} \g_2 \big)^{\prime}
=I_{m_2}-\I A E_0^{-1}  \big(\wt E^{-1} \g_2 \big)^{\prime},
\end{align} 
respectively, we see that $E_0^{-1}$ also maps differentiable functions with a square-summable derivative
into differentiable functions with a square-summable derivative.
Thus, the same is valid for $E=\wt E E_0$ and for $E^{-1}$.
\epr
\begin{Rk}\label{RkDop1}
Formulas $\wt E^{-1}A=A\wt E^{-1}$ and \eqref{y13} for $\wt E^{-1}$ and formulas above for
$E_0^{-1}$ yield a useful equality
\begin{align} &      \label{r1}
\big(E^{-1}\g_1\big)(0)=\big(E_0^{-1}\wt E^{-1}\g_1\big)(0)=\g_1(0)=0.
\end{align} 
\end{Rk}
Now, we construct a representation of the fundamental solution $w$
of the system
\begin{align} &      \label{i9}
\frac{d}{dx}w(x,z)=\I zj\g(x)^*\g(x)w(x,z), \quad w(0,z)=I_m.
\end{align} 
For that purpose we introduce operators
\begin{align} &      \label{i10}
S:=E^{-1}\big(E^*\big)^{-1}, \quad \Pi:= \begin{bmatrix}
\Phi_1 & \Phi_2
\end{bmatrix}, \quad 
\Phi_k \in B\big(\mathbb{C}^{m_k}, \, L^2_{m_2}(0, \, l)\big);
\\
&      \label{i10'}
\big(\Phi_1 f\big)(x)=\Phi_1(x)f, \quad
\Phi_1(x):=\big(E^{-1}\g_1\big)(x);  \quad  \Phi_2 f=I_{m_2}f\equiv f;
\end{align}
where $E$ is constructed (for the given $\g$) in Proposition \ref{PnSimN} and $\g_1$
is the left $m_2\times m_1$ block of $\g$.
We also introduce the transfer matrix function in Lev Sakhnovich form \cite{SaL1, SaL2, SaL3}
\begin{align} &      \label{i11}
w_A(z):=I_m+\I zj\Pi^*S^{-1}(I-zA)^{-1}\Pi.
\end{align} 
We shall need  the reductions of the operators above (and the matrix function $w_A$ corresponding 
to those reductions):
\begin{align} &      \label{i12}
\big(P_{\xi}f\big)(x)=f(x) \quad (0<x<\xi), \quad P_{\xi}\in B\Big(L^2_{m_2}(0, \, \eT), \, L^2_{m_2}(0, \, \xi)\Big),
\\ &      \label{i13}
 A_{\xi}:=P_{\xi}AP_{\xi}^*, \quad S_{\xi}:=P_{\xi}SP_{\xi}^*,
\\ &      \label{i14}
w_A(\xi,z):=I_m+\I zj\Pi^*P_{\xi}^*S_{\xi}^{-1}(I-zA_{\xi})^{-1}P_{\xi}\Pi, \quad 0<{\xi}\leq \eT.
\end{align}
\begin{Tm}\label{FundSol}
Let $\g$ be determined by \eqref{3.1}, where $u$ is the fundamental solution
of the Dirac system \eqref{1.1}  with a square-summable potential $v$. 
 Then, the fundamental solution $w$ given by  \eqref{i9} 
 admits representation
 \begin{align} &      \label{i15}
w(\xi,z)=w_A(\xi,z),
\end{align} 
where $w_A(\xi,z)$ is defined by  \eqref{i14}.
\end{Tm}
\bpr 
Formulas \eqref{i3}, \eqref{i10} and \eqref{i10'} imply that
\begin{align} &      \label{i21}
 \Pi f=(E^{-1}\g)f.
\end{align}  
It is immediate from the definition \eqref{3.4} of $K$ that
 \begin{align} &      \label{i16}
K^*=-\I \int_x^{\eT}\g(x)j\g(t)^* \cdot dt, \quad K-K^*=\I \g(x) j\int_0^{\eT}\g(t)^*\cdot dt.
\end{align} 
According to Proposition \ref{PnSimN} we have $K=EAE^{-1}$. Since $K=EAE^{-1}$,
taking into account  \eqref{i10} and \eqref{i21}, we rewrite
the second equality in \eqref{i16} in the form of the operator identity
 \begin{align} &      \label{i17}
AS-SA^*=\I\Pi j \Pi^*.
\end{align} 
Hence, we may use the Method of Operator Identities \cite{SaL1, SaL2, SaL3}.
We need now to  show the applicability of
the Continuous Factorization Theorem (see \cite[p. 40]{SaL3})
or, more conveniently, its corollary  \cite[Theorem 1.20]{SaSaR}.
Completely similar to the cases in \cite{SaSaR} we see that conditions $(i)$ and $(ii)$ of    
\cite[Theorem 1.20]{SaSaR} are satisfied. It remains only to derive that
$\Pi^*P_{\xi}S_{\xi}^{-1}P_{\xi}\Pi$ is absolutely continuous (i.e., condition  $(iii)$ of    
\cite[Theorem 1.20]{SaSaR} holds)
and that 
 \begin{align} &      \label{i18}
\left(\Pi^*P_{\xi}^*S_{\xi}^{-1}P_{\xi}\Pi\right)^{\prime}=H(\xi)=\g(\xi)\g(\xi)^*,
\end{align} 
in order to prove that $w_A$ satisfies the differential system in \eqref{i9}.

Since the operator $E$ is invertible, triangular, and has Hilbert-Schmidt kernel, we see that $E^{-1}$ is also triangular.
Taking into account that $E^{\pm 1}$ are lower triangular operators, we obtain
\begin{align} &      \label{i19}
P_{\xi}EP_{\xi}^*P_{\xi}=P_{\xi}E, \qquad  \big(E^{-1}\big)^*P_{\xi}^*=P_{\xi}^*P_{\xi}\big(E^{-1}\big)^*P_{\xi}^*.
\end{align}  
The first equality in \eqref{i19} yields $P_{\xi}EP_{\xi}^*P_{\xi}E^{-1}P_{\xi}^*=P_{\xi}P_{\xi}^*$,
that is, 
$$P_{\xi}E^{-1}P_{\xi}^*=(P_{\xi}EP_{\xi}^*)^{-1}.$$  
Hence, formulas  (\ref{i10}),
 (\ref{i13}), and \eqref{i19} lead us to
\begin{align} &      \label{i20}
 S_{\xi}^{-1}=E_{\xi}^*E_{\xi}, \qquad E_{\xi}:=P_{\xi}EP_{\xi}^*.
\end{align}  
Finally, from \eqref{i21}, \eqref{i19}, and \eqref{i20} we derive that
\begin{align} &      \label{i22}
\Pi^*P_{\xi}^*S_{\xi}^{-1}P_{\xi}\Pi=\int_0^{\xi}\g(\zeta)\g(\zeta)^*d\zeta
\end{align}  
(i.e., $\Pi^*P_{\xi}^*S_{\xi}^{-1}P_{\xi}\Pi$ is absolutely continuous and \eqref{i18} is valid).
Hence, $w_{A}$ satisfies the system in \eqref{i9} and, furthermore, the normalization
\begin{align} &      \label{i22'}
\lim_{x\to 0}w_A(x,z)=I_m
\end{align}  
easily follows from \eqref{i14} and \eqref{i20}.
\epr

\noindent Since \eqref{i17} holds we say that the triple $\{A, S,  \Pi\}$ forms an $S$-node~\cite{SaL1, SaL2, SaL3}.
\begin{Cy}\label{FSD} Let $u(x,z)$ be the fundamental solution   of a Dirac system with the square-summable potential
$v$ and let $\g$ be given by  \eqref{3.1}. Then $u(x,z)$ admits representation
\begin{align} &      \label{i23}
u(x,z)=\E^{\I xz}u(x,0)w_{A}(x,2z).
\end{align} 
 Here $w_A$ has the form  \eqref{i14},  where the $S$-node $\{A, \, S, \, \Pi\}$, which determines $w_A$,    
is given in \eqref{i1}, \eqref{i10}, and \eqref{i10'}.
\end{Cy}
\bpr
According to \eqref{1.1} and Theorem \ref{FundSol} we have
\begin{align}       \nn
\big(\E^{\I xz}u(x,0)w_{A}(x,2z)\big)^{\prime}=&(\I z I_m+\I jV(x)+2 \I z u(x,0)j\g(x)^*\g(x)u(x,0)^{-1})
\\ & \label{i24} \,\, \times
\E^{\I xz}u(x,0)w_{A}(x,2z).
\end{align} 
Writing $u(x,0)$ in the block form and taking into account \eqref{i5}, we derive
\begin{align}       &\label{i25}
u(x,0)=\begin{bmatrix} \b(x)\\ \g(x)\end{bmatrix}, \quad u(x,0)j\g(x)^*=\begin{bmatrix} 0 \\ -I_{m_2} \end{bmatrix}.
\end{align}
From \eqref{i4} we obtain $u(x,0)^{-1}=ju(x,0)^*j$. Thus, in view of \eqref{i24} and \eqref{i25} we see that
\begin{align}       \nn
\big(\E^{\I xz}u(x,0)w_{A}(x,2z)\big)^{\prime}=&\left(\I z I_m+\I jV(x)-2\I z \begin{bmatrix} 0 &0 \\ 0 & I_{m_2} \end{bmatrix} \right)\E^{\I xz}u(x,0)
\\  \label{i26}  \times w_{A}(x,2z)
=&(\I zj+\I j V(x))\E^{\I xz}u(x,0)w_{A}(x,2z).
\end{align} 
Relations \eqref{i22'} and \eqref{i26} yield \eqref{i23}.
\epr
\section{Solution of the inverse problem}
\setcounter{equation}{0}
Here,  we may follow the lines of \cite[Sections 3 and 4]{FKRS3} without any essential changes.
The high-energy asymptotics of $\vp$ is given by the following theorem.
\begin{Tm}\label{TmHea} Assume that $\vp \in \cln(\eT,z)$ and the potential $v$ of the corresponding   Dirac system  \eqref{1.1} 
is square-summable on $(0, \, \eT)$. Then $($uniformly with respect to $\Re(z))$ we have
\begin{align} \label{i27}&
\vp(z)=2\I z\int_0^{\eT}\E^{2\I xz}\Phi_1(x)dx+O\left(2z \E^{2\I \eT z}/ \sqrt{\Im(z)}\right), \quad \Im (z) \to \infty.
\end{align} 
\end{Tm}
\bpr To prove the theorem, we consider the matrix function
\begin{align} &\label{i28}
\clu(z)=\begin{bmatrix}
I_{m_1} & \vp(z)^*
\end{bmatrix}\big(j-
w_A({\eT},2z)^*jw_A({\eT},2z)\big)
\begin{bmatrix}
I_{m_1} \\ \vp(z)
\end{bmatrix}.
\end{align} 
It easily follows from  (\ref{i14}) and   (\ref{i17}) (see, e.g., \cite[p. 24]{SaSaR})
that
\begin{align} \label{i33}&
w_A(\eT, z)^*jw_A(\eT,z)=j+\I(z - \ov z)\Pi^*(I-\ov z A^*)^{-1}S^{-1}(I-zA)^{-1}\Pi ,
\end{align} 
and so we derive $\clu(z)\geq 0$. Because of  (\ref{i4}), (\ref{i23}), and  (\ref{i28}) we have
\begin{align} &\label{i29}
\clu(z)=I_{m_1}-\vp(z)^*\vp(z)- \E^{\I \eT(\ov z-  z)}\begin{bmatrix}
I_{m_1} & \vp(z)^*
\end{bmatrix}
u({\eT},z)^*ju({\eT},z)
\begin{bmatrix}
I_{m_1} \\ \vp(z)
\end{bmatrix}.
\end{align} 
We note that \eqref{2.2} yields
\begin{align} &\label{i30}
\begin{bmatrix}
I_{m_1} \\ \vp(z)
\end{bmatrix}=u(\eT,z)^{-1}
\clp(z)
\big(\begin{bmatrix}
I_{m_1} & 0
\end{bmatrix}u(\eT,z)^{-1}\clp(z)\big)^{-1}.
\end{align} 
Taking into account   (\ref{i30}), we rewrite \eqref{i29} as
\begin{align} \nn
\clu(z)=&I_{m_1}-\vp(z)^*\vp(z)-\E^{\I \eT(\ov z-  z)}
\Big(\big(\begin{bmatrix}
I_{m_1} & 0
\end{bmatrix}u(\eT,z)^{-1}\clp(z)\big)^{-1}\Big)^*
 \\ & \,\, \label{i31}
\times
\clp(z)^*j\clp(z)
\big(\begin{bmatrix}
I_{m_1} & 0
\end{bmatrix}u(\eT,z)^{-1}\clp(z)\big)^{-1}.
\end{align} 
Recall that $\clu(z)\geq 0$. Hence, from  (\ref{2.1}) and  (\ref{i31}) we see that
\begin{align} &\label{i32}
0\leq \clu(z)\leq I_{m_1}, \quad \vp(z)^*\vp(z)\leq I_{m_1}.
\end{align} 
Now, formulas  (\ref{i28}),  (\ref{i33}), and  (\ref{i32}) imply that
\begin{align} &\label{i34}
 2\I (\ov z-z) \begin{bmatrix}
I_{m_1} & \vp(z)^*
\end{bmatrix}
\Pi^*(I-2\ov z A^*)^{-1}S^{-1}(I-2zA)^{-1}\Pi 
\begin{bmatrix}
I_{m_1} \\ \vp(z)
\end{bmatrix}\leq I_{m_1}.
\end{align} 
Since $S$ is positive and boundedly invertible, inequality  (\ref{i34}) yields
\begin{align} &\label{i35}
\left\| (I-2zA)^{-1}\Pi 
\begin{bmatrix}
I_{m_1} \\ \vp(z)
\end{bmatrix}\right\| \leq C/\sqrt{\Im z} \quad {\mathrm{for}}\,\, 
{\mathrm{some}} \quad C>0.
\end{align} 
After applying $-\I\Phi_2^*$ to the operator on the
left-hand side  of  (\ref{i35}), we derive
\begin{align} &\label{i36}
-\I\Phi_2^* (I-2zA)^{-1}\Phi_2\vp(z)=\I\Phi_2^* (I-2zA)^{-1}\Phi_1+
O\left(\frac{1}{\sqrt{\Im(z)}}\right).
\end{align} 
Using \eqref{y5'} we see that
\begin{align} \label{i37}&
\Phi_2^*(I-2zA)^{-1}f=\int_0^{\eT}\E^{2\I(x-\eT)z}f(x)dx, 
\\  \label{i37'}&
  \Phi_2^*(I-2zA)^{-1}\Phi_2=\frac{\I}{2z}\big(\E^{-2\I \eT z}-1\big)I_{m_2}.
\end{align} 
Because of  (\ref{i36})--(\ref{i37'}), we have
\begin{align} \label{i38}&
\frac{1}{2z}\big(\E^{-2\I \eT z}-1\big)\vp(z)=\I \E^{-2\I \eT z}\int_0^{\eT}\E^{2\I x z}\Phi_1(x)dx+
O\left(\frac{1}{\sqrt{\Im(z)}}\right).
\end{align} 
Since $\vp$ is non-expansive, we see  from  (\ref{i38}) that  (\ref{i27}) holds.
\epr
\begin{Cy}\label{cyHea} Let $\vp$ be the Weyl function of Dirac system  \eqref{1.1}
on  $[0, \, \infty)$, where the potential $v$    
is locally square-summable. Then  we have
\begin{align} \label{Repr}&
\vp(z)=2\I z\int_0^{\infty}\E^{2\I xz}\Phi_1(x)dx, \quad \Im (z) >0.
\end{align} 
\end{Cy}
\bpr Since $\vp$ is analytic and non-expansive in $\mathbb{C}_+$, 
for any $\ve>0$ it admits (see, e.g., \cite[Theorem V]{WP} or a slightly more convenient for us reformulation
\cite[Theorem E.11]{SaSaR}) a representation
\begin{align} \label{repr0}&
\vp(z)=2\I z\int_0^{\infty}\E^{2\I xz}\Phi (x)dx, \quad \Im (z) >\ve>0,
\end{align} 
where $\E^{-2\ve x}\Phi(x) \in L^2_{m_2\times m_1}(0, \, \infty)$. Because of  (\ref{i27})
and  (\ref{repr0}) we obtain
\begin{align}\nn
\psi(z):&=\int_0^{\eT}\E^{2\I(x-\eT)z}\big(\Phi_1 (x)-\Phi(x)\big)dx
\\  \label{4.15} &
=\int_{\eT}^{\infty}\E^{2\I(x-\eT )z}\Phi(x)dx+O\big(1/\sqrt{\Im (z)}\big).
\end{align} 
From  (\ref{4.15}) we see that $\psi(z)$ is bounded in some half-plane
$\Im(z)\geq \eta_0>0$. Clearly, $\psi(z)$ is bounded also in the half-plane
$\Im(z)< \eta_0$. Since $\psi$ is analytic and bounded in $\mathbb{C}$
and  tends to zero on some rays, we have
\begin{align} \label{4.16} &
\psi(z)=\int_0^{\eT}\E^{2\I(x-\eT)z}\big(\Phi_1 (x)-\Phi(x)\big)dx \equiv 0.
\end{align} 
It follows from  (\ref{4.16}) that $\Phi_1(x)\equiv \Phi(x)$ on all 
finite intervals $[0, \, \eT]$. Hence,  (\ref{repr0}) implies  (\ref{Repr}).
\epr

\begin{Rk}\label{RkL}
According to the proof of Corollary \ref{cyHea}, we have $\Phi_1\equiv \Phi$, and so $ \Phi_1(x)$ does not depend
on $\eT$ for $\eT>x$. 
Furthermore,  the proof of  Corollary \ref{cyHea}
implies also that
$\E^{-\ve x}\Phi_1(x) \in L^2_{m_2\times m_1}(0, \, \infty)$ for any $\ve >0$.
\end{Rk}

Using representation \eqref{Repr}, we uniquely recover $v$ from $\vp$.
Indeed, taking into account Plancherel Theorem and Remark \ref{RkL},
we apply inverse Fourier transform to formula  (\ref{Repr}) and derive
\begin{align} \label{5.1}&
\Phi_1\Big(\frac{x}{2}\Big)=\frac{1}{\pi}\E^{x\eta}{\mathrm{l.i.m.}}_{a \to \infty}
\int_{-a}^a\E^{-\I x\xi}\frac{\vp(\xi+\I \eta)}{2\I(\xi +\I \eta)}d\xi, \quad \eta >0.
\end{align} 
Here l.i.m. stands for the entrywise limit in the norm of  $L^2(0,b)$, 
$\, 0<b \leq \infty$. 
(Note that if we put additionally $\Phi_1(x)=0$ for $x<0$, equality  (\ref{5.1})
holds for l.i.m. as the entrywise limit in $L^2(-b,b)$.)
Thus, for any fixed interval $(0, \, \eT)$ the corresponding operators $S$ and $\Pi$ are recovered from $\vp$. 

Since the Hamiltonian $H$
is recovered from $S$ and $\Pi$ via formula  (\ref{i18}), and $H=\g^*\g$,
we recover also $\g$. 
First, for that purpose,  we recover  the so called Schur coefficient:                                                                                        
\begin{align} \label{5.2}&
\left(\begin{bmatrix}
0 &I_{m_2}
\end{bmatrix}H\begin{bmatrix}
0 \\ I_{m_2}
\end{bmatrix}\right)^{-1}
\begin{bmatrix}
0 &I_{m_2}
\end{bmatrix}H\begin{bmatrix}
I_{m_1} \\ 0
\end{bmatrix}=(\g_2^*\g_2)^{-1}\g_2^*\g_1=\g_2^{-1}\g_1.
\end{align} 
Here we used the inequality $\det \g_2\not=0$, which follows from the second identity
in \eqref{i5}. The second identity
in \eqref{i5} yields also
$$I_{m_2}-(\g_2^{-1}\g_1)(\g_2^{-1}\g_1)^*=\g_2^{-1}(\g_2^{-1})^*,$$
which implies that the left-hand side of this equality is invertible.
Taking into account  $\det \g_2\not=0$,  we rewrite $\g_1$
 in the form $\, \g_1=\g_2 (\g_2^{-1}\g_1)$
and the identity 
\eqref{i6} in the form $\g_2^{\prime}=\g_1^{\prime}(\g_2^{-1}\g_1)^*$. Therefore, we obtain 
\begin{align} &      \nn
\g_2^{\prime}=(\g_2 (\g_2^{-1}\g_1))^{\prime}(\g_2^{-1}\g_1)^*, \quad {\mathrm{i.e.,}}\\
 &      \label{p18}
\g_2^{\prime}=
\g_2(\g_2^{-1}\g_1)^{\prime}(\g_2^{-1}\g_1)^*
\big(I_{m_2}-(\g_2^{-1}\g_1)(\g_2^{-1}\g_1)^*\big)^{-1},
\end{align} 
and recover $\g_2$ from \eqref{p18} and the initial condition $\g_2(0)=I_{m_2}$.
Finally, we recover $\g_1$ from $\g_2$ and $\g_2^{-1}\g_1$.

In order to recover $\b$ from $\g$, we partition $\b$ into two blocks $\b=\begin{bmatrix}
\b_1 & \b_2
\end{bmatrix}$, where $\b_k$ ($k=1,2$) is an $m_1\times m_k$ matrix function.
We put
\begin{align} \label{5.3}&
\wt \b=\begin{bmatrix}
I_{m_1} & \g_1^*(\g_2^*)^{-1}
\end{bmatrix}.
\end{align} 
Because of  (\ref{i5}) and  (\ref{5.3}), we have $\b j \g^*=\wt \b j \g^*=0$, and so
\begin{align} \label{5.4}&
\b(x)= \b_1(x)\wt \b(x).
\end{align} 
It follows from  (\ref{1.1}) and  (\ref{3.1}) that 
\begin{align} \label{5.4'}&
\b^{\prime}(x)=\I v(x)\g(x),
\end{align} 
which implies
\begin{align} \label{5.5}&
\b^{\prime}j\b^*=0, \qquad \b^{\prime}j\g^*=-\I v.
\end{align} 
Formula  (\ref{5.4}) and the first relation in  (\ref{i5})  lead us to
\begin{align} \label{5.6}&
 \wt \b j \wt \b^*=\b_1^{-1}(\b_1^*)^{-1}.
\end{align} 
From \eqref{5.4} we also derive that
\begin{align} \nn&
\b^{\prime}j\b^*=\b_1^{\prime}(\wt \b j\wt \b^*)\b_1^*+\b_1(\wt \b^{\prime}j\wt \b^*)\b_1^*.
\end{align}
Taking into account the first relation in  (\ref{5.5}) and formula \eqref{5.6}, we rewrite the
relation above:
\begin{align} \label{5.6'}&
\b_1^{\prime}\b_1^{-1}+\b_1(\wt \b^{\prime}j\wt \b^*)\b_1^*=0.
\end{align}
According to  (\ref{1.3}), \eqref{5.6},  and  (\ref{5.6'}),  $\b_1$ satisfies the first order differential equation
(and initial condition):
\begin{align} \label{5.7}&
\b_1^{\prime}=-\b_1(\wt \b^{\prime}j\wt \b^*)( \wt \b j \wt \b^*)^{-1}, \quad \b_1(0)=I_{m_1}.
\end{align} 
Thus, $\b_1$ and $\b$ are successively recovered from $\g$.
The potential $v$ is recovered from $\b$ and $\g$ via the second equality in \eqref{5.5}.
In this way, we recover $v$ on any interval $[0, \, \eT]$, therefore, on the whole semiaxis.
We proved the following theorem.
 \begin{Tm} \label{TmIP} Let $\vp$ be the Weyl function of Dirac system  \eqref{1.1} on $[0, \, \infty)$,
 where the potential $v$ is  locally square-summable.
 Then $v$ can be uniquely recovered from $\vp$ via the formula
\begin{align} &      \label{5.8}
v(x)=\I \b^{\prime}(x)j\g(x)^*.
\end{align}  
Here $\b$ is recovered from $\g$ using  \eqref{5.3},  \eqref{5.4} and  \eqref{5.7}$;$ $\g$ is recovered
from the Hamiltonian $H$ using  \eqref{5.2} and   \eqref{p18}$;$ the Hamiltonian is given by
 \eqref{i18}, $\Pi$  from  \eqref{i18} is expressed via $\Phi_1(x)$ in formula
 \eqref{i10'}, and  $S$ is the unique solution of \eqref{i17}. Finally, $\Phi_1(x)$ is recovered from $\vp$ using  \eqref{5.1}.
 \end{Tm}
 \begin{Rk}\label{RkSxi} It follows from \eqref{i13} and \eqref{i17} that the operator identities
 \begin{align} &      \label{r2}
A_{\xi}S_{\xi}-S_{\xi}A_{\xi}^*=\I P_{\xi}\Pi j (P_{\xi}\Pi)^*, \quad 0<\xi \leq \eT,
\end{align} 
where $A$ is given in \eqref{i1}, $A_{\xi}= P_{\xi}AP_{\xi}^*$, and $\Pi$ is given by \eqref{i10'}, hold. 
The uniqueness of the operators $S_{\xi}$ satisfying these identities is proved on p.311 in \cite{SaSaR}.
Moreover, it is easy to see that the proof of \cite[Proposition~3.2]{FKS1} works also for the case, where
$\psi$ and $\wt \psi$ are differentiable functions with the square-summable derivatives.
Thus, recalling \eqref{r1} and formulas (3.16) and (3.17) in \cite[Proposition 3.2]{FKS1}, we see
that $S_{\xi}$ given by
\begin{align} &      \label{r3}
S_{\xi}=I-\frac{1}{2}\int_0^{\xi}\int_{|x-t|}^{x+t}\Phi_1^{\prime}\left(\frac{\zeta+x-t}{2}\right)
\Phi_1^{\prime}\left(\frac{\zeta+t-x}{2}\right)^*d\zeta \, \cdot \, dt
\end{align} 
satisfies \eqref{r2}. Hence, $S_{\xi}$ of the form \eqref{r3} is the unique solution of \eqref{r2},
and we may recover $S_{\xi}$ $($considered in Theorem \ref{TmIP}$)$ from $\Phi_1$ in this way.
 \end{Rk}
Using Theorem \ref{TmIP} we modify Borg-Marchenko-type Theorem 2.52 from \cite{SaSaR} for the case
of the locally square-summable potentials. We note that seminal publications by F. Gesztesy and B. Simon \cite{GS, GeSi, Si0}
gave rise to a series of interesting results on the high energy asymptotics of the Weyl functions and local
Borg-Marchenko-type uniqueness theorems. Recall that the high energy asymptotics of the Weyl functions
is given (for our case) in Theorem \ref{TmHea}.
 
\begin{Tm}\label{BM} Let $\vp$ and $\wh \vp$ be  Weyl functions
of twoDirac systems on $[0, \, \eT]$ $($or on $[0, \, \infty))$ with  square-summable $($locally square-summable$)$ potentials,
which are denoted by $\, v\,$ and $\,\wh v\, $, respectively. Suppose that on some ray
$\Re z$~$=$~$c \Im z$, where $\,c \in \BR$ and $\Im z>0$, the equality
\begin{align}     \label{sar1}&
\|\vp(z)-\wh \vp(z)\|=O(\E^{2\I \zeta z }) \quad (\Im z \to \infty) 
\end{align} 
holds  for
all $0<\zeta < l \,\, ( l<\eT<\infty)$.
Then we have
\begin{align}     \label{sar2}&
v(x)=\wh v(x), \qquad 0<x< l.
\end{align} 
\end{Tm}
\bpr Since Weyl functions are non-expansive, it is immediate that
the inequality
\begin{align}     \label{sar3}&
\|\E^{-2\I \zeta z }\big(\vp(z)-\wh \vp(z)\big)\|\leq c_1 \E^{2\zeta |z| }, \quad \Im z \geq c_2>0
\end{align} 
is valid for some $c_1$ and $c_2$. It is apparent also that
the matrix function $\E^{-2\I \zeta z }\big(\vp(z)-\wh \vp(z)\big)$ is bounded
on  the line $\Im z =c_2$. Furthermore, formula  (\ref{sar1}) implies that
$\E^{-2\I \zeta z }\big(\vp(z)-\wh \vp(z)\big)$ is bounded on the ray $\Re z=c \Im z$.
Therefore, applying the Phragmen-Lindel\"of theorem (e.g., its version  \cite[Corollary E.7]{SaSaR}) in the angles
generated by the line $\Im z =c_2$ and the ray $\Re z=c \Im z$ ($\Im z \geq c_2$),
we see that
\begin{align}     \label{sar4}&
\|\E^{-2\I \zeta z }\big(\vp(z)-\wh \vp(z)\big)\|\leq c_3, \quad \Im z \geq c_2>0.
\end{align} 
Let functions associated with $\wh \vp$ be written with a hat (e.g.,
$\wh v, \,\wh \Phi_1$).  Because of formula  (\ref{i27}), its analog for $\wh \vp$,
$\wh \Phi_1$ and the inequality  (\ref{sar4}), we have
\begin{align}     \label{sar5}&
\left\| \int_0^{\zeta} \E^{2\I  (x- \zeta)z}\big(\Phi_1(x)-\wh \Phi_1(x)\big)dx   \right\|
\leq c_4, \quad \Im z \geq c_2>0.
\end{align} 
Clearly, the left-hand side of  (\ref{sar5}) is bounded in the half-plane
$\Im z<c_2$ and tends to zero on some rays. Thus, we derive
\begin{align}     \label{sar6}&
\int_0^{\zeta} \E^{2\I  (x- \zeta)z}\big(\Phi_1(x)-\wh \Phi_1(x)\big)dx   \equiv 0,
\quad {\mathrm{i.e.}}, \quad  \Phi_1(x)\equiv \wh \Phi_1(x) \quad
(0<x<\zeta).
\end{align} 
Since  (\ref{sar6}) holds for all $\zeta< l$, we obtain $ \Phi_1(x)\equiv \wh \Phi_1(x)$ for $0<x<  l$. In view of Theorem \ref{TmIP}, the last identity implies  (\ref{sar2}).
\epr

{\bf Acknowledgement.}
The research was supported by the
Austrian Science Fund (FWF) under Grant  No. P24301.
The author is grateful to F.~Gesztesy for his question, which initiated this paper.

\end{document}